\documentclass[a4paper,11pt]{amsart}

\usepackage{times}
\usepackage{amsthm,latexsym,amssymb,amsmath,amsrefs}
\usepackage[bookmarks=true]{hyperref}   % Hyperref in DVI and PDF
\usepackage{listings}
\usepackage{nicefrac}
\usepackage{ccicons}

\hypersetup{
    colorlinks,
    linkcolor={red!50!black},
    citecolor={blue!50!black},
    urlcolor={blue!80!black}
}

\newcommand{\R}{\mathbb{R}}
\newcommand{\C}{\mathbb{C}}
\newcommand{\Q}{\mathbb{Q}}
\newcommand{\Z}{\mathbb{Z}}
\newcommand{\<}{\!\left\langle}
\renewcommand{\>}{\right\rangle}
\newcommand{\tr}{\mathrm{tr}}
\newcommand{\pf}{\mathrm{pf}}
\newcommand{\sign}{\mathrm{sign}}
\newcommand{\adj}{\mathrm{adj}}
\newcommand{\Padj}{\mathrm{Padj}}
\newcommand{\Mat}{\mathrm{Mat}}
\newcommand{\Skew}{\mathrm{Skew}}

\renewcommand{\O}{\mathrm{O}}

\renewcommand{\phi}{\varphi}

\newtheorem{cor}{Corollary}
\newtheorem{prop}{Proposition}
\newtheorem{lemma}{Lemma}

\theoremstyle{definition}
\newtheorem{ex}{Example}

\newtheorem{definition}{Definition}

\usepackage{xcolor}

\definecolor{codegreen}{rgb}{0,0.6,0}
\definecolor{codegray}{rgb}{0.5,0.5,0.5}
\definecolor{codepurple}{rgb}{0.58,0,0.82}
\definecolor{backcolour}{rgb}{0.95,0.95,0.92}

\lstdefinestyle{mystyle}{
    backgroundcolor=\color{backcolour},   
    commentstyle=\color{codegreen},
    keywordstyle=\color{magenta},
    numberstyle=\tiny\color{codegray},
    stringstyle=\color{codepurple},
    basicstyle=\ttfamily\footnotesize,
    breakatwhitespace=false,         
    breaklines=true,                 
    captionpos=b,                    
    keepspaces=true,                 
    numbers=none,                    
    numbersep=5pt,                  
    showspaces=false,                
    showstringspaces=false,
    showtabs=false,                  
    tabsize=2
}

\title{The Faddeev-LeVerrier algorithm and the Pfaffian}
\author{Christian B\"ar}
\date{\today}
\address{Christian B\"ar, 
Institut f\"ur Mathematik,
Universit\"at Potsdam, 
D-14476, Potsdam, Germany
}
\urladdr{\href{https://www.math.uni-potsdam.de/baer}{https://www.math.uni-potsdam.de/baer}}
\email{\href{mailto:cbaer@uni-potsdam.de}{cbaer@uni-potsdam.de}}
\thanks{\ccCopy\, 2021. This manuscript version is made available under the CC-BY-NC-ND 4.0 license \url{https://creativecommons.org/licenses/by-nc-nd/4.0/}}

\begin{document}

\begin{abstract}
We adapt the Faddeev-LeVerrier algorithm for the computation of characteristic polynomials to the computation of the Pfaffian of a skew-symmetric matrix.
This yields a very simple, easy to implement and parallelize algorithm of computational cost $\O(n^{\beta+1})$ where $n$ is the size of the matrix and $\O(n^{\beta})$ is the cost of multiplying $n\times n$-matrices, $\beta\in[2,2.37286)$.
We compare its performance to that of other algorithms and show how it can be used to compute the Euler form of a Riemannian manifold using computer algebra.
\end{abstract}

\maketitle

\section{Introduction}

The computation of the determinant of a large $n\times n$-matrix $A$ is a challenging task because naively using the formula $\det(A) = \sum_\sigma \sign(\sigma) a_{1\sigma(i)}\cdots a_{n\sigma(n)}$ would require $\O((n-1)n!)$ many multiplications since the are $n!$ permutations $\sigma$ of $\{1,\ldots,n\}$.
This is practical only for very small matrices.

The Gauss algorithm allows to transform a matrix in triangular form with $\O(n^3)$ many operations.
The determinant can then be read off.
However, the Gauss algorithm includes divisions by matrix entries, which requires a pivoting strategy to avoid numerical instabilities.
More seriously, it cannot be carried out if the matrix takes entries in a commutative ring in which we cannot divide.

It is therefore desirable to find a fast division-free algorithm.
The classical Faddeev-LeVerrier algorithm provides precisely this; its computational complexity is of the order $\O(n^{\beta+1})$ and it avoids any division by a matrix entry.
Here $\O(n^\beta)$ is the computational cost of the multiplication of two $n\times n$-matrices.
The precise value of $\beta$ is unknown but it does lie in the interval $[2,2.37286)$, compare the discussion at the beginning of Section~\ref{sec:performance}.

Moreover, the Faddeev-Leverrier algorithm is very simple and easy to implement.
One essentially has to carry out $n$ matrix multiplications; hence using any software which parallelizes matrix multiplication will automatically parallelize the whole algorithm.
In Section~\ref{sec:FL} we will recall the algorithm and present the probably shortest possible derivation.
It even simplifies the approach presented in \cite{Hou1998}.
However, no claim of originality is made for this part.
See \cite{Barnett, Givens, Gower, HWV, WangLin} for various extensions of the algorithm.

The Faddeev-LeVerrier algorithm is division free in the sense that no divisions by matrix entries are required. 
But we do have to divide by integers.
For this reason we allow matrices with entries in general commutative $\Q$-algebras $R$.
We could as well allow any torsion-free commutative ring $R$ because such an $R$ embeds into its rationalization $\Q\otimes_\Z R$ via $r\mapsto 1\otimes r$.
So the algorithm can be used if $R=\Z$, $R=\Q$, $R=\R$, $R=\C$, or if $R$ is any field of characteristic $0$.
If $R$ has torsion, e.g.\ if $R$ is a finite field, then different algorithms are required, see e.g.\ \cite{Berkowitz}.

Our interest lies in the fast computation of the Pfaffian.
The Pfaffian is defined for skew-symmetric matrices and is a polynomial in its entries which squares to the determinant.
The Pfaffian may be less prominent than the determinant but it has important applications in physics, in combinatorics, and in geometry.
As to physics, see the introduction of \cite{GBRB2011,Wimmer2012} and the references therein and the application to the topological charge of a disordered nanowire in \cite{Wimmer2012}.
In combinatorics, the Pfaffian of the directed adjacency matrix of a suitable graph yields the number of perfect matchings.
In Section~\ref{sec:application} we will see an application to differential geometry.

As for determinants, the computation of a Pfaffian using the definition to be given in \eqref{eq:defPfaff} would be way too slow.
There are established algorithms of order $\O(n^3)$ which transform the matrix in a normal form from which the Pfaffian can be read off.
This is similar to computing the determinant using Gauss elimination.
These algorithms are well suited for the numerical treatment of real or complex matrices and, again, cannot be applied to matrices with entries in commutative rings.

We adapt the Faddeev-LeVerrier algorithm to compute the Pfaffian in Section~\ref{sec:FLPfaff} and obtain an algorithm of order $\O(n^{\beta+1})$.
Again, it is very simple, easy to implement and can easily be parallelized.
The derivation of this algorithm requires some theoretical material about Pfaffians which we provide in Section~\ref{sec:Pfaff} since it seems hard to find a good presentation in the literature.
The algorithm is again division free in the sense that no divisions by matrix entries are required. 

In Section~\ref{sec:performance} we compare the performance of our algorithm to that of other established algorithms.
Based on an earlier version of this paper, our algorithm has been implemented in SageMath~9.3 which allows for an easy comparison with the previous implementation using the definition of the Pfaffian.

The main application we have in mind is to compute the Pfaffian of matrices taking entries in algebras of symbolic expressions as typically used in computer algebra.
In Section~\ref{sec:application} we show how the algorithm can be used to compute the Euler form of a Riemannian manifold using SageMath.
In this case, the entries of the matrix are formal differential forms of mixed even degree.
Thus we are dealing with a commutative $\Q$-algebra in which we cannot divide.
In particular, algorithms based on transforming the matrix in a canonical form do not apply.
Moreover, multiplication in this algebra is very costly as it involves expensive simplification routines.
Our algorithm appears to provide the best known approach to deal with this kind of scenario.

\subsection*{Acknowledgment}
The author wants to thank Zachary Hamaker for pointing out the complexity results for matrix multiplication and Darij Grinberg and an unknown referee for many comments and suggestions on how to improve the presentation.
Moreover, he is grateful for financial support by SPP 2026 funded by Deutsche Forschungsgemeinschaft.

\section{Faddeev-LeVerrier algorithm for the characteristic polynomial}
\label{sec:FL}

Let $R$ be a commutative $\Q$-algebra.
In particular, this covers the case $R=\Q$, $R=\R$, and $R=\C$.
In Section~\ref{sec:application} we will consider a more sophisticated algebra.
Let $A\in\Mat(n,R)$ where $\Mat(n,R)$ denotes the set of all $n\times n$-matrices with entries in $R$.
We want to compute the coefficients of its characteristic polynomial $\chi(t)=\det(tI-A)$, in particular its determinant.

\subsection{Derivation of the algorithm}
To start, we write $\chi(t)=\sum_{j=0}^n c_{n-j}t^j$. 
Then $c_0=1$, $c_1=-\tr(A)$, and $c_n=(-1)^n\det(A)$.
Recall that $\det(A)\cdot I=A\cdot\adj(A)$ where $\adj(A)$ is the adjugate matrix of $A$.
The entries of $\adj(A)$ are given by determinants of $(n-1)\times(n-1)$-submatrices of $A$ and hence are polynomials of degree $n-1$ in the entries of $A$.
Applying this to the matrix $tI-A$, we get
\begin{equation*}
\chi(t)\cdot I = (tI-A)\sum_{j=0}^{n-1} t^j N_{n-j},
%\label{eq:resolvent}
\end{equation*}
where the $N_{n-j}$ are $n\times n$-matrices with entries in $R$.
Comparing coefficients in
\begin{align*}
\sum_{j=0}^n c_{n-j}t^j \, I
% =
% \chi(t) \, I
&=
(tI-A)\sum_{j=0}^{n-1} t^j N_{n-j} \\
&=
-AN_n + \sum_{j=1}^{n-1} t^j (N_{n-j+1}-AN_{n-j}) + t^nN_1
\end{align*}
yields
\begin{align*}
c_0\cdot I &= N_1, \\
c_k\cdot I &= N_{k+1}-AN_k\quad\mbox{ for }k=1,\ldots,n-1, \\
c_n\cdot I &= -AN_n
\end{align*}
and hence
\begin{align}
N_1 &= I, \label{eq:N1}\\
N_{k+1} &= AN_k + c_k\,I\quad\mbox{ for }k=1,\ldots,n-1, \label{eq:Nk}\\
0 &= AN_n + c_n\cdot I.\label{eq:Nn}
\end{align}
If we knew the coefficients $c_k$, equations \eqref{eq:N1} and \eqref{eq:Nk} would provide a recursive procedure to determine the matrices $N_k$.
Equation~\eqref{eq:Nn} would serve as an additional check.

Up to this point the derivation of the algorithm is the usual one to be found in many places in the literature.
To get more information about the $c_k$, one traditionally uses Newton's formulas for the power sums and the elementary symmetric polynomials (see e.g.\ \cite{Gantmacher_1}*{Ch.~IV, \S~5} or \cite{Householder}*{Sec.~6.7}).
Hou \cite{Hou1998} applies the Laplace transform to the matrix exponential.

Instead, we just use Jacobi's formula which states that the logarithmic derivative of the determinant is given by the trace of the logarithmic derivative of the matrix.
Denote by $\tfrac{d}{dt}$ the formal derivative of polynomials in $R[t]$.
Using \eqref{eq:Nk} and \eqref{eq:N1} we find
\begin{align*}
\dot\chi(t)
&=
\frac{d}{dt}\det(tI-A)
=
\tr\big(\tfrac{d}{dt}(tI-A)\cdot\adj(tI-A)\big) \\
&=
\tr(\adj(tI-A))
=
\sum_{j=0}^{n-1} t^j \tr(N_{n-j}) \\
&=
\sum_{j=0}^{n-2} t^j \tr(AN_{n-j-1}+c_{n-j-1}I) + t^{n-1}\tr(N_1) \\
&=
\sum_{j=0}^{n-2} t^j \big(\tr(AN_{n-j-1})+nc_{n-j-1}\big) + nt^{n-1}.
\end{align*}
On the other hand,
\begin{equation*}
\dot\chi(t)
=
\sum_{j=0}^{n-1} (j+1)c_{n-j-1}t^j.
\end{equation*}
Comparing coefficients again, we obtain
\begin{equation*}
\tr(AN_k) + nc_k = (n-k)c_k
\end{equation*}
and hence 
\begin{equation}
c_k = -\tfrac{1}{k}\tr(AN_k)
\end{equation}
for $k=1,\ldots,n$.
Inserting this into \eqref{eq:Nk} gives us the recursion procedure
\begin{align}
N_1 &= I, \label{eq:NN1}\\
N_{k+1} &= AN_k -\tfrac{1}{k}\tr(AN_k)\,I\quad\mbox{ for }k=1,\ldots,n \label{eq:NNk}
\end{align}
and $N_{n+1}=0$ is the additional check.

\subsection{Implementation}

The following Python code implements the algorithm.
After an initialization like
\medskip
\begin{lstlisting}[language=Python, caption={Initialization}]
import numpy as np
n = 50                                 # matrix size
I = np.identity(n)
A = np.random.randn(n,n).astype('i8')  # replace by your matrix
\end{lstlisting}
and running the algorithm by
\medskip
\begin{lstlisting}[language=Python, caption={Faddeev-LeVerrier algorithm}]
c = [1,-A.trace()]
N = A+c[1]*I
for k in range(2,n+1):
    M = np.matmul(A,N)
    c.append(-M.trace()/k)
    N = M + c[k]*I 
\end{lstlisting}
the list ${c}$ contains the coefficients of the characteristic polynomial, ${c}=[c_0,c_1,\ldots,c_n]$.
The matrix ${N}$ is then the matrix $N_{n+1}=0$.
If one skips the computation of ${N}$ for the last iteration with $k=n$, then $N$ will end up to be $N_n$, which by \eqref{eq:Nn} coincides with $(-1)^{n+1}\adj(A)$.
Thus the algorithm computes the adjugate matrix of $A$ (and hence the inverse matrix if $A$ is invertible) along the way.

\section{The Pfaffian}
\label{sec:Pfaff}

\subsection{Definition of the Pfaffian}
Denote the space of skew-symmetric $n\times n$-matrices by $\Skew(n,R) := \{A\in\Mat(n,R) \mid A^\top=-A\}$.
Let $n=2m$ be even and let $A\in\Skew(n,R)$.
The Pfaffian of $A=(a_{ij})$ is defined as 
\begin{equation}
\pf(A) = \sum_P \sign(P)\cdot a_{i_1j_1}\cdots a_{i_mj_m}
\label{eq:defPfaff}
\end{equation}
where the sum is taken over all perfect matchings $P$ of $\{1,\ldots,n\}$.
Here a perfect matching is a partition of the set $\{1,\ldots,n\}$ into a disjoint union of $m$ sets with two elements each, $\{1,\ldots,n\}=\{i_1,j_1\}\sqcup\cdots\sqcup\{i_m,j_m\}$.
If we use the convention that each pair is ordered by $i_k<j_k$ then the sign of $P$ is defined as the sign of the permutation $\begin{bmatrix}1&2&\cdots&n-1&n\\ i_1&j_1&\cdots&i_m&j_m\end{bmatrix}$.
Note that this is well defined because changing the order in which the pairs occur does not change the sign of the corresponding permutation.

\begin{ex}
Let
$$
J=
\begin{pmatrix}
\, \boxed{\begin{matrix} 0&  1  \\ -1&  0\end{matrix}}&  &  \\ 
 &  \ddots  &  \\ 
 &  &  \boxed{\begin{matrix} 0&  1  \\ -1&  0\end{matrix}} 
\end{pmatrix}
.
$$
In this case only the perfect matching $\{1,2\}\sqcup\cdots\sqcup\{n-1,n\}$ contributes.
Hence $\pf(J)=1$.
\end{ex}

\subsection{Properties of the Pfaffian}
If $A$ is skew-symmetric then 
\begin{equation}
\pf(A)^2 = \det(A)
\label{eq:pfdet3}
\end{equation}
and if $B$ is any $n\times n$-matrix then 
\begin{equation}
\pf(B^\top AB) = \det(B)\cdot \pf(A),
\label{eq:pftrafo}
\end{equation}
see e.g.\ \cite{Artin}*{Thms.~3.27 and 3.28}.
Hence the Pfaffian is a polynomial in the entries of a skew-symmetric matrix which squares to the determinant.

If $A$ is invertible in $\Mat(n,R)$ then $\det(A)$ is invertible in $R$ and hence $\pf(A)$ is invertible in $R$ too, with inverse $\frac{\pf(A)}{\det(A)}$.

If we interchange two rows and the corresponding columns in $A$, then $\pf(A)$ changes sign.
This is a consequence of \eqref{eq:pftrafo} since this change is obtained by choosing $B$ to be a transposition matrix.
For example, if we want to interchange the first and second rows and columns we apply \eqref{eq:pftrafo} with
$$
B=
\begin{pmatrix}
\, \boxed{\begin{matrix} 0&  1  \\ 1&  0\end{matrix}}&  &  & \\
 & 1 & & & \\
 & &  \ddots  &  \\ 
 & & & 1
\end{pmatrix}
.
$$
In particular, if $A$ has two identical rows then $\pf(A)=-\pf(A)$ and hence $\pf(A)=0$.

\subsection{The Laplace expansion}
We need an analogue to the Laplace expansion of a determinant.
For $A\in\Skew(n,R)$ and $i\neq j\in\{1,\ldots,n\}$ denote by $A\< i,j\>\in\Skew(n-2,R)$ the matrix obtained from $A$ by removing the $i$-th and $j$-th row and column.
%This is a skew-symmetric $(n-2)\times(n-2)$-matrix.

\begin{prop}\label{prop:Laplace}
Fix $i\in\{1,\ldots,n\}$.
Then
$$
\pf(A) = \sum_{j< i} (-1)^{i+j} \cdot a_{ij} \cdot \pf(A\<i,j\>) + \sum_{j> i} (-1)^{i+j+1} \cdot a_{ij} \cdot \pf(A\<i,j\>).
$$
\end{prop}

\begin{proof}
Any perfect matching $P$ of $\{1,\ldots,n\}$ can be uniquely written in the form $\{1,\ldots,n\}=\{i,j\}\sqcup P'$ where $j\neq i$ and $P'=\{i_2,j_2\}\sqcup\cdots\sqcup\{i_m,j_m\}$ is a perfect matching of $\{1,\ldots,n\}\setminus \{i,j\}$.
Then $\sign(P)=(-1)^{i+j+1}\cdot\sign(P')$. 
Hence
\begin{align*}
\pf(A) 
&=
\sum_{j> i}\sum_{P'} (-1)^{i+j+1} \sign(P') a_{ij}\cdot a_{i_2j_2}\cdots a_{i_mj_m} \\
&\quad +
\sum_{j< i}\sum_{P'} (-1)^{i+j+1} \sign(P') a_{ji}\cdot a_{i_2j_2}\cdots a_{i_mj_m} \\
&=
\sum_{j> i} (-1)^{i+j+1} \cdot a_{ij} \cdot \pf(A\<i,j\>) 
+
\sum_{j< i} (-1)^{i+j} \cdot a_{ij} \cdot \pf(A\<i,j\>) .
\qedhere
\end{align*}
\end{proof}

\begin{ex}
Choosing $i=2$ this yields for a general skew-symmetric $4\times 4$-matrix:
\begin{align*}
\pf&
\begin{pmatrix}
0         & a_{12}         & a_{13}         & a_{14} \\
-a_{12}   & 0              & a_{23}         & a_{24} \\
-a_{13}   & -a_{23}        & 0              & a_{34} \\
-a_{14}   & -a_{24}        & -a_{34}        & 0
\end{pmatrix} \\
&=
-(-a_{12})\cdot\pf\begin{pmatrix}0 & a_{34} \\ -a_{34} & 0\end{pmatrix}
+ a_{23}\cdot\pf\begin{pmatrix}0 & a_{14} \\ -a_{14} & 0\end{pmatrix}
- a_{24}\cdot\pf\begin{pmatrix}0 & a_{13} \\ -a_{13} & 0\end{pmatrix} \\
&=
a_{12}a_{34} + a_{23}a_{14} - a_{24}a_{13}.
\end{align*}
\end{ex}

\subsection{The Pfaff-adjugate matrix}

\begin{definition}
We call the matrix $\Padj(A)=(b_{ij})\in\Mat(n,R)$ the \emph{Pfaff-adjugate matrix} of $A$ where
$$
b_{ij} =
\begin{cases}
0, & \mbox{ if }i=j,\\
(-1)^{i+j}  \cdot \pf(A\<i,j\>), & \mbox{ if } i<j,\\
(-1)^{i+j+1}  \cdot \pf(A\<i,j\>), & \mbox{ if } i>j.
\end{cases}
$$
\end{definition}

Note that $\Padj(A)$ is skew-symmetric as well.
The terminology is justified by

\begin{cor}\label{cor:adjugate}
Let $A\in\Skew(n,R)$. 
Then
$$
A\cdot \Padj(A) = \pf(A)\cdot I.
$$
In particular, if $A$ is invertible then $A^{-1}=\frac{1}{\pf(A)}\Padj(A)$.
\end{cor}

\begin{proof}
We compute the entries of $A\cdot \Padj(A)$.
For the diagonal entries we get
\begin{align*}
\sum_{j=1}^n a_{ij}b_{ji}
&=
\sum_{j<i} a_{ij}\cdot (-1)^{i+j}  \cdot \pf(A\<i,j\>) + \sum_{j>i} a_{ij}\cdot (-1)^{i+j+1}  \cdot \pf(A\<i,j\>) \\
&=
\pf(A)
\end{align*}
by Proposition~\ref{prop:Laplace}.

Now let $i\neq k$.
We only consider the case $i<k$, the case $i>k$ being analogous.
Let $\tilde A$ be the matrix obtained from $A$ by replacing the $k$-th row and column by the $i$-th row and column, respectively.
Then $\tilde A$ is again a skew-symmetric $n\times n$-matrix but since it has two identical rows we have $\pf(\tilde A)=0$.
Denote the entries of $\tilde A$ by $\tilde a_{\mu\nu}$.
Notice $ a_{ij}=\tilde a_{ij}=\tilde a_{kj}$ for all $j$.
We compute
\begin{align*}
\sum_{j=1}^n a_{ij}b_{jk}
&=
\sum_{j<k} a_{ij}\cdot (-1)^{k+j}  \cdot \pf(A\<k,j\>) + \sum_{j>k} a_{ij}\cdot (-1)^{k+j+1}  \cdot \pf(A\<k,j\>) \\
% &=
% \sum_{j<k} \tilde a_{ij}\cdot (-1)^{k+j}  \cdot \pf(\tilde A\<k,j\>) + \sum_{j>k} \tilde a_{ij}\cdot (-1)^{k+j+1}  \cdot \pf(\tilde A\<k,j\>) \\
&=
\sum_{j<k} \tilde a_{kj}\cdot (-1)^{k+j}  \cdot \pf(\tilde A\<k,j\>) + \sum_{j>k} \tilde a_{kj}\cdot (-1)^{k+j+1}  \cdot \pf(\tilde A\<k,j\>) \\
% &\stackrel{(*)}{=}
% \sum_{j<i} \tilde a_{ij}\cdot (-1)^{k+j}  \cdot (-1)^{k-i-1}  \cdot \pf(\tilde A\<i,j\>) \\
% &\quad
% + \sum_{i<j< k} \tilde a_{ij}\cdot (-1)^{k+j} \cdot (-1)^{k-i-2}  \cdot \pf(\tilde A\<i,j\>) \\
% &\quad
% + \sum_{j>k} \tilde a_{ij}\cdot (-1)^{k+j+1}  \cdot (-1)^{k-i-1}   \cdot \pf(\tilde A\<i,j\>) \\
% &=
% \sum_{j<i} \tilde a_{ij}\cdot (-1)^{i+j+1}  \cdot \pf(\tilde A\<i,j\>) 
% + \sum_{j>i} \tilde a_{ij}\cdot (-1)^{i+j}  \cdot \pf(\tilde A\<i,j\>) \\
&=
-\pf(\tilde A)
= 0
\end{align*}
again by Proposition~\ref{prop:Laplace}.
% The marked equation holds because $\tilde A\<i,j\>$ is obtained from $\tilde A\<k,j\>$ by interchanging $k-i-1$ rows and columns if $j$ is smaller than $i$ or larger than $k$ and by interchanging $k-i-2$ rows and columns if $j$ lies between $i$ and $k$.
\end{proof}

\begin{cor}
Let $A\in\Skew(n,R)$. 
Then
$$
\pf(A)\cdot\Padj(A) = \adj(A).
$$
\end{cor}

\begin{proof}
From
$$
A\cdot\adj(A)
=
\det(A)\cdot I
=
\pf(A)^2\cdot I
= 
\pf(A)\cdot A\cdot\Padj(A)
$$
we conclude
$$
A\cdot(\adj(A)-\pf(A)\Padj(A)) = 0.
$$
The entries of $E(A):=\adj(A)-\pf(A)\Padj(A)$ are universal polynomials of degree at most $n-1$ in the entries of $A$.
Note that $E(A)$ vanishes whenever $A$ is invertible.
If $R=\R$ then $E(A)$ vanishes for all $A$ because invertible matrices are dense in $\Skew(n,\R)$.
Thus $E=0$ as a polynomial in the entries of $A$.
Hence the assertion follows for general $R$.
\end{proof}

We will need the Pfaffian version of Jacobi's formula.

\begin{lemma}\label{lem:JacobiPfaff}
Let $A(t)$ be a skew-symmetric $n\times n$-matrix with entries in $R[t]$.
Then
$$
\tfrac{d}{dt}\pf(A(t)) = \tfrac12\tr(\dot A(t)\cdot \Padj(A(t))).
$$
\end{lemma}

\begin{proof}
We compute
\begin{align*}
\pf(A(t))\cdot\tr\big(\dot A(t)\cdot \Padj(A(t))\big)
&=
\tr\big(\dot A(t)\cdot \adj(A(t))\big) \\
&=
\tfrac{d}{dt}\det(A(t)) \\
&=
2\,\pf(A(t))\cdot\tfrac{d}{dt}\pf(A(t))
\end{align*}
and hence
$$
\pf(A(t))\cdot\Big(\tfrac{d}{dt}\pf(A(t)) - \tfrac12\tr\big(\dot A(t)\cdot \Padj(A(t))\big)\Big) = 0.
$$
The expression $E(t) := \tfrac{d}{dt}\pf(A(t)) - \tfrac12\tr\big(\dot A(t)\cdot \Padj(A(t))\big)$ is a polynomial in~$t$,
$$
E(t) = \sum_k p_k t^k.
$$
The coefficients $p_k$ are universal polynomials (depending only on $n$ but not on $R$) in the coefficients of the entries of $A(t)\in\Mat(n,R[t])$,
$$
p_k = p_k(a_{ij\ell};\quad i,j=1,\ldots,n,\, \ell = 0,\ldots,k+1)
$$
where $A(t)=(\sum_\ell a_{ij\ell}t^\ell )$.
These universal polynomials have rational coefficients.

If $R=\R$ and $A(t_0)$ is invertible for some $t_0\in\R$ then $A(t)$ is invertible for $t$ near $t_0$ and hence $E(t)=0$ for $t$ near $t_0$.
Thus $E=0$ as a polynomial.
Since invertible matrices are dense in $\Skew(n,\R)$ the lemma follows if $R=\R$.

This implies that the universal polynomials $p_k$ vanish as polynomials and hence the lemma holds for arbitrary $R$.
\end{proof}

Darij Grinberg pointed out that one can also prove Lemma~\ref{lem:JacobiPfaff} without referring to determinants by a direct computation similar to the proof of Jacobi's formula.

\section{The Faddeev-LeVerrier algorithm for the Pfaffian}
\label{sec:FLPfaff}

Now we are ready to adapt the Faddeev-LeVerrier algorithm to compute the Pfaffian.

\subsection{Derivation of the algorithm}
We consider the \emph{Pfaffian characteristic polynomial}\footnote{Compare also the \emph{quasi-characteristic polynomial} in \cite{Krivo}.} of our skew-symmetric $n\times n$-matrix $A$ where $n=2m$,
$$
\Psi(t) = \pf(tJ+A) .
$$
This is a polynomial in $t$ of degree $m$.
We expand it $\Psi(t)=\sum_{j=0}^m c_{m-j}t^j$. 
Then $c_0=1$ and we are interested in computing $c_m=\Psi(0)=\pf(A)$.

The entries of $\Padj(A)$ are given by Pfaffians of $(n-2)\times(n-2)$-submatrices of $A$ and hence are polynomials of degree $m-1$ in the entries of $A$.
Applying Corollary~\ref{cor:adjugate} to the matrix $tJ+A$, we get
\begin{equation}
\Psi(t)\cdot I = (tJ+A)\cdot\sum_{j=0}^{m-1} t^j N_{m-j},
\label{eq:resolventPf}
\end{equation}
where $N_{m-j}\in\Mat(n,R)$.
Comparing coefficients in
\begin{align*}
\sum_{j=0}^m c_{m-j}t^j \, I
=
\Psi(t) \, I
&=
(tJ+A)\sum_{j=0}^{m-1} t^j N_{m-j} \\
&=
AN_m + \sum_{j=1}^{m-1} t^j (JN_{m-j+1}+AN_{m-j}) + t^mJN_1
\end{align*}
yields
\begin{align*}
c_0\cdot I &= JN_1, \\
c_k\cdot I &= JN_{k+1}+AN_k\quad\mbox{ for }k=1,\ldots,m-1, \\
c_m\cdot I &= AN_m
\end{align*}
and hence
\begin{align}
N_1 &= -J, \label{eq:N1Pf}\\
JN_{k+1} &= -AN_k + c_k\,I\quad\mbox{ for }k=1,\ldots,n-1, \label{eq:NkPf}\\
0 &= AN_m - c_m\cdot I.\label{eq:NnPf}
\end{align}
Lemma~\ref{lem:JacobiPfaff} yields
\begin{align*}
\dot\Psi(t)
&=
\tfrac{d}{dt}\pf(tJ+A)
=
\tfrac12\tr(\tfrac{d}{dt}(tJ+A)\cdot\Padj(tJ+A)) \\
&=
\tfrac12\tr(J\Padj(tJ+A))
=
\frac{1}{2}\sum_{j=0}^{m-1} t^j \tr(JN_{m-j})
\end{align*}
and hence, by \eqref{eq:NkPf} and \eqref{eq:N1Pf},
\begin{align*}
\dot\Psi(t)
&=
% \frac12\sum_{j=0}^{m-1} t^j \tr(JN_{m-j})
% =
\frac12\sum_{j=0}^{m-2} t^j \tr(-AN_{m-j-1}+c_{m-j-1}I) + \tfrac12 t^{m-1}\tr(I) \\
&=
\frac12\sum_{j=0}^{m-2} t^j \big(-\tr(AN_{m-j-1})+nc_{m-j-1}\big) + mt^{m-1}.
\end{align*}
Comparing coefficients with
\begin{equation*}
\dot\Psi(t)
=
\sum_{j=0}^{m-1} (j+1)c_{m-j-1}t^j
\end{equation*}
yields
\begin{equation*}
-\tfrac12\tr(AN_k) + mc_k = (m-k)c_k
\end{equation*}
and hence 
\begin{equation}
c_k = \tfrac{1}{2k}\tr(AN_k)
\end{equation}
for $k=1,\ldots,n$.
Inserting this into \eqref{eq:NkPf} gives us the recursion procedure
\begin{align}
N_1 &= -J, \label{eq:NN1Pf}\\
N_{k+1} &= JAN_k -\tfrac{1}{2k}\tr(AN_k)\,J\quad\mbox{ for }k=1,\ldots,n \label{eq:NNkPf}.
\end{align}

\subsection{Implementation}

The following Python code implements the Faddeev-LeVerrier algorithm for the Pfaffian.
After initializing

\medskip
\begin{lstlisting}[language=Python, caption={Initialization for the computation of the Pfaffian}]
import numpy as np
m = 20                    # matrix size
n = 2*m
J = np.zeros((n,n)).astype('i8')
for k in range(m):
    J[2*k,2*k+1] = 1
    J[2*k+1,2*k] = -1
A = np.random.randn(n,n).astype('i8')
A = A - A.transpose()     # make skewsymmetric
\end{lstlisting}
and running the algorithm by
\medskip
\begin{lstlisting}[language=Python, caption={Faddeev-LeVerrier algorithm for the Pfaffian}, label=list:FLPfaff]
c = 1
N = -J
for k in range(1,m+1):
    M = np.matmul(A,N)
    c = M.trace()/(2*k)
    if k<m: N = np.matmul(J,M) - c*J
\end{lstlisting}
the variable $c$ contains $\pf(A)$.

\subsection{Remarks}
After termination of the algorithm the variable $N$ contains $N_m$, for which we have by \eqref{eq:NnPf}
$$
\pf(A)\cdot I = A\cdot N_m
$$
Indeed, $N_m$ is the Pfaff-adjugate matrix of $A$.
So again, the algorithm computes the inverse of an invertible skew-symmetric matrix along the way.

The auxiliary matrix $J$ occurring in the Pfaffian characteristic polynomial replaces the identity matrix $I$ in the usual characteristic polynomial.
Unlike $I$, the matrix $J$ is not invariant under similarity transformations and therefore the Pfaffian characteristic polynomial is less canonical than the usual characteristic polynomial.
This is the reason why the other coefficients $c_k$ for $k<m$ in the Pfaffian characteristic polynomial are of little interest.

In fact, the algorithm works with any skew-symmetric matrix $J$ satisfying $J^2=-I$.
Another natural choice would be
$$
J = 
\begin{pmatrix}
0 & I \\
-I & 0
\end{pmatrix}.
$$
One just has to replace $c_0=1$ by $c_0=\pf(J)$ and \eqref{eq:NN1Pf} gets replaced by $N_1=-\pf(J)J$.

\section{Performance}
\label{sec:performance}

There are $(n-1)!!$ many perfect matchings of $\{1,\ldots,n=2m\}$.
Thus a direct implementation of the Pfaffian based on \eqref{eq:defPfaff} would require $\O(n(n-1)!!)$ many multiplications which is unpractical unless the matrix is very small.
 
Obviously, the multiplication of two $n\times n$-matrices can be performed at computational cost $\O(n^3)$.
Surprisingly, Strassen in \cite{Strassen} found a way to perform matrix multiplication at cost $\O(n^{\beta})$ with $\beta=\log_27$.
The optimal exponent $\beta$ is unknown to date.
It clearly satisfies $\beta\ge2$ and after several improvements by various authors (see e.g.\ \cite{LeGall,Schoenhage}) Alman and Vassilevska Williams showed $\beta<2.37286$, see \cite{Alman}.

Since in our algorithm there is just one loop with essentially one matrix multiplication at each iteration, we find that our algorithm is of order $\O(n^{\beta+1})$, hence better than $\O(n^{3.37286})$.

\subsection{Performance comparison with other algorithms}
It seems that the current versions Maple~2021.0 and Mathematica~12.3 do not have the Pfaffian implemented.
SageMath has had an implementation based directly on \eqref{eq:defPfaff} for many years.
Based on an earlier version of this paper, our algorithm has been implemented in version 9.3 of SageMath.
If applicable, i.e.\ if the matrix take entries in a torsionfree ring, our version of the Faddeev-LeVerrier algorithm performs much better.

The following tests were carried out with SageMath~9.3 on a Linux system with an AMD Ryzen 9 3900 processor and 128 GB RAM.
Table~\ref{tab:Sage} shows the average wall time in seconds for the computation of the Pfaffian of a skew-symmetric $n\times n$-matrix with rational random entries, once computed with our algorithm and once with the original implementation.

\begin{table}[h!]
\begin{tabular}{|c|c|c|}
\hline
$n$ & our algorithm & original \\
\hline\hline
10 & 0.00012 & 0.025 \\
\hline
16 &  0.00041 & 68\\
\hline
20 & 0.00085 & $\infty$ \\
\hline
\end{tabular}
\vspace{2mm}
\caption{Comparison of algorithms for the method \texttt{pfaffian} in SageMath~9.3}
\label{tab:Sage}
\end{table}

In this test all computations were exact with rational numbers represented by fractions.
For numerical computation, more efficient algorithms are available.
Wimmer approaches the problem in \cite{Wimmer2012} by first transforming the matrix into a partial skew-symmetric tridiagonal form from which the Pfaffian can be easily determined.
This is similar in spirit to using the Gauss algorithm to compute the determinant.
Wimmer implemented two methods to carry out the transformation, the Householder and the Parlett-Reid algorithm.

Table~\ref{tab:Wimmer} shows the averaged wall times in seconds in the Python implementation of Wimmer's and our algorithm when applied to skew-symmetric $n\times n$-numpy-matrices with random float64-entries.

\begin{table}[h!]
\begin{tabular}{|c|c|c|c|}
\hline
$n$ & our algorithm & Wimmer with Householder & Wimmer with Parlett-Reid \\
\hline\hline
10 &0.00017 & 0.00026  & 0.00016\\
\hline
100 & 0.0059 & 0.0032 &  0.0014\\
\hline
500 & 2.2 & 0.10 & 0.050\\
\hline
\end{tabular}
\vspace{2mm}
\caption{Comparison with Wimmer's algorithm}
\label{tab:Wimmer}
\end{table}

For not too large matrices ($n\le100$) our version of the Faddeev-LeVerrier algorithm is comparable to Wimmer's algorithm but for larger matrices Wimmer is faster. 
This is to be expected as the computational cost of Wimmer's algorithm is of the order $\O(n^3)$.
According to Table~II in \cite{Wimmer2012}, the Fortran implementation of Wimmer's algorithm clearly beats other existing similar algorithms such as the one in LAPACK and the one by Gonz\'alez-Ballestero et al. \cite{GBRB2011}.

For even larger matrices such as $n=1000$ neither ours nor Wimmer's algorithm returns a result because of numerical inaccuracies, at least in the Python implementation.

A combinatorial algorithm based on graph-theoretic methods has been proposed in \cite{MSV2004}, see also \cite{Rote2001}.
It would be interesting to compare its performance to that of the Faddeev-LeVerrier algorithm.

\subsection{Optimizations}
The costly part in Listing~\ref{list:FLPfaff} is the matrix multiplication $M=A*N$ which has to be carried out $m=\tfrac{n}{2}$ times.
Matrix multiplication is ideally suited for parallelization as the entries can be computed independently from each other.
Thus using software in which matrix multiplication is parallelized will automatically parallelize the Faddeev-LeVerrier algorithm.

Note that the matrix multiplication $J*M$ is cheap because it just amounts to swapping rows in $M$ and changing signs in half of them.

\section{An application in differential geometry}
\label{sec:application}

In this last section we want to illustrate the usefulness of the algorithm for symbolic computation.
Our version of the Faddeev-LeVerrier algorithm is division free in the sense that we never have to divide by an entry of the skew-symmetric matrix $A$ whose Pfaffian we are computing. 
This not only avoids numerical instabilities but also allows us to consider matrices with entries in general commutative $\Q$-algebras in which division may not be possible.
We will make use of this now.

Let $M$ be an oriented compact differentiable manifold of even dimension \mbox{$n=2m$.}
Let $g$ be a Riemannian metric on $M$, i.e.\ $g$ provides each tangent space $T_pM$ with a scalar product which depends smoothly on the base point $p\in M$.
The associated curvature tensor $R\colon T_pM\times T_pM\times T_pM\to T_pM$ has the property that $(X,Y,Z,W) \mapsto g(R(X,Y)Z,W)$ is skew-symmetric in $X$ and $Y$ and also in $Z$ and $W$. 

Fix $p\in M$ and let $s=(s_1,\ldots,s_n)$ be a positively oriented basis of $T_pM$.
Then 
$$
\Omega_s:=(g(R(\cdot,\cdot)s_i,s_j))
$$ 
is a skew-symmetric matrix with entries in $\Lambda^2 T_p^*M$, the space of skew-symmetric $2$-forms on $T_pM$.
We consider the commutative algebra 
$$
R 
%= \Lambda^\mathrm{even}T_p^*M 
:= \bigoplus_{k=0}^m \Lambda^{2k} T_p^*M
$$ 
where multiplication in this algebra is given by the wedge product.
The Pfaffian of $\Omega_s\in\Mat(n,R)$ is a homogeneous polynomial of degree $m$ in the entries of $\Omega_s$ and hence a form of degree $2m=n$, i.e.\ $\pf(\Omega_s)\in\Lambda^nT_p^*M$.

The Pfaffian $\pf(\Omega_s)$ depends on the choice of $s$ however.
Let $e=(e_1,\ldots,e_n)$ be another positively oriented basis of $T_pM$ and let $B=(b_{ij})$ be the transformation matrix characterized by $s_i=\sum_k b_{ki}e_k$.
Then $\Omega_s = B^\top \cdot \Omega_e \cdot B$ and hence by \eqref{eq:pftrafo} 
\begin{equation}
\pf(\Omega_s)= \det(B)\pf(\Omega_e).
\label{eq:EulerTrafo}
\end{equation}
If $s$ and $e$ are both orthonormal bases, then $B\in\mathrm{SO}(n)$ and thus $\pf(\Omega_s)= \pf(\Omega_e)$.
We call $\chi(R) := \pf(\Omega_e)$ the \emph{Euler form} if $e$ is orthonormal.
This yields a well-defined $n$-form on $M$ whose importance comes from the Gauss-Bonnet-Chern theorem (\cite{Chern1944}):
$$
\int_M \chi(R) = (2\pi)^m\cdot \mbox{Euler-Poincar\'e characteristic}(M).
$$
If $e$ is orthonormal then we easily find for the symmetric matrix $G_s := (g(s_i,s_j))$ that
$$
G_s = B^\top\cdot B
$$
and therefore $\det(B) = \sqrt{\det(G_s)}$.
By \eqref{eq:EulerTrafo} we can write the Euler form for a general basis $s$ as
$$
\chi(R) = \det(B)^{-1}\cdot\pf(\Omega_s) = \det(G_s)^{-\nicefrac12} \cdot \pf(\Omega_s).
$$
To compute the Euler form in examples one can use the SageMath code developed in \cite{Jung2020}.
Multiplications in the algebra $R=\bigoplus_{k=0}^m \Lambda^{2k} T_p^*M$ are computationally expensive.
So the classical Faddeev-LeVerrier algorithm and our modification of it are a good choice to compute $\det(G_s)$ and $\pf(\Omega_s)$, respectively.

This approach can be easily generalized to compute the Euler form of an arbitrary vector bundle (rather than the tangent bundle $TM$) equipped with a semi-Riemannian metric, compare \cite{Chern1963}.

\subsection*{Note added in proof}
Fredrik Johansson pointed out that the Faddeev-Leverrier algorithm for the determinant has been modified by Preparata and Sarwate \cite{PS} in such a way that it uses only $\O(\sqrt{n})$-many matrix multiplications.
This improves its computational cost to $\O(n^{\beta+0.5} + n^3)$ multiplications, see \cite{J} for more details.
The modification can also be applied to our version of the algorithm computing the Pfaffian.

\end{document}